\documentclass[12pt]{amsart}
\usepackage{epsf, amssymb, xy, latexsym, amsmath, graphics}
\newtheorem{theorem}{Theorem}
\newtheorem{corollary}{Corollary}
\newtheorem{lemma}{Lemma}
\newtheorem{proposition}{Proposition}
\def \1{^{-1}}
\def \2{^{-2}}
\def \3{^{-3}}

\def\t{\tau}
\def\r{\rho}
\def\s{\sigma}
\def \pr{^\prime}
\newcommand{\C}{{\mathbb  C}} 
\newcommand{\N}{{\mathbb  N}} 
\newcommand{\Z}{{\mathbb  Z}} 
 
\newcommand{\B}{{\mathbb B}} 
\newcommand{\cf}{{\mathcal F}}
\newcommand{\cg}{{\mathcal G}} 
\newcommand{\ch}{{\mathcal H}}
\newcommand{\ck}{{\mathcal K}}
\newcommand{\cl}{{\mathcal L}} 
\newcommand{\cm}{{\mathcal M}} 
 
\newcommand{\co}{{\mathcal O}} 

\newcommand{\cp}{{\mathcal P}}

\newcommand{\cv}{{{\mathcal V}an}}
\renewcommand{\b}{{\beta}}
\renewcommand{\c}{{\gamma}}
\newcommand{\abdots}{{a,b_0,\dots,b_m}}
\newcommand{\bdots}{{b_0,\dots,b_m}}
\newcommand{\braid}{{\mathbf B}}
\newcommand{\pure}{{\mathbf P}}

\newcommand{\pzero}{{\mathbb P}} 
\newcommand{\pone}{{\pzero^1}} 
\newcommand{\ptwo}{{\pzero^2}}

\newcommand{\pn}{{\pzero^n}}

\newcommand{\cn}{{\C^n}}

\newcommand{\orbfg}{\pi_1^{orb}}
\newcommand{\ok}{\to}
\newcommand{\ise}{\Rightarrow}
\newcommand{\onto}{\twoheadrightarrow}
\newcommand{\eksi}{\backslash}
\newcommand{\ssol}{\langle\!\langle}
\newcommand{\ssag}{\rangle\!\rangle}
\newcommand{\comb}[2]{\left(\!\!{\begin{array}{c}#1\\#2 \end{array}}\!\!\right)}
\newcommand{\comp}{{\,\mbox{\small o}\,}}
\newcommand{\col}[2]{{\begin{array}{l}{#1}\\{#2}\end{array}}}
\title{On Branched Galois Coverings $(\B_1)^n\ok \pn$} 
\author{A.Muhammed Uluda{\u g}} 
\begin{document} 
\begin{abstract}  
For any $n>1$, we construct examples branched Galois coverings $M\ok\pn$ where $M$ is one of 
$(\pone)^n$, $\cn$  and $(\B_1)^n$, and  $\B_1$ is the 1-ball.  
In terms of orbifolds, this amounts to giving examples of orbifolds over $\pn$ uniformized by $M$.
\end{abstract} 
\maketitle  

\noindent
\textbf{1. Introduction.}
In contrast with the considerable  literature on the orbifolds over 
$\ptwo$ uniformized by the 2-ball $\B_2$
(see~\cite{yoshida},~\cite{deligne2},~\cite{hirzebruch} and references therein), 
not much is known about which orbifolds over $\pn$ are uniformized by other symmetric spaces. 
In this article, we apply a simple orbifold-covering technique 
to construct some orbifolds over the projective space $\pn$ uniformized either by 
$(\pone)^n$, $\cn$ or $(\B_1)^n$. Our main result is the following theorem.
\begin{theorem}\label{main}
Let $(n,b)$ be a pair of coprime integers with $n\geq 2$.
There exists a Galois covering  
$\bigl(D^{(b)}_{n,1}\bigr)^n\ok \pn$ of degree $n!b^{n^2-n}$
branched  along an irreducible degree-$2b(n-1)$ 
hypersurface $D_n^{(b)}\subset\pn$ where
$D^{(b)}_{n,1}\subset D_n^{(b)}$ is a curve of  euler number $e=b^{n-1}(n+1+b-nb)$.
 \end{theorem}
\noindent
For $b=1$, the hypersurface $D_n^{(1)}$ is the discriminant hypersurface, and 
$D_{n,1}^{(1)}\simeq \pone$ is a rational normal curve. 
In this case one obtains the well-known branched Galois covering $(\pone)^n\ok\pn$.
The subvarieties $D_n^{(b)}$ and $D_{n,1}^{(b)}$ are the liftings respectively of $D_n^{(1)}$ and
$D_{n,1}^{(1)}$ by an abelian branched self-covering 
$[Z_0,\dots,Z_n]\in\pn \ok [Z_0^b,\dots, Z_n^b]\in \pn$.
For $(n,b)\in \{(3,2),(2,3)\}$ one has $e\bigl(D_{n,1}^{(b)}\bigr)=0$, and the 
universal covering of $\bigl(D^{(b)}_{n,1}\bigr)^n $ is $\C^n$. 
The curve $D_2^{(3)}=D_{2,1}^{(3)} $ is the nine-cuspidal sextic. 
For $b>1$ and $(n,b)\notin \{(3,2),(2,3)\}$ one has  $e\bigl(D_{n,1}^{(b)}\bigr)<0$, and the 
universal covering of $\bigl(D^{(b)}_{n,1}\bigr)^n $ is $(\B_1)^n$.
\par
In case $(n,b)=(2,3)$, the claim of Theorem~\ref{main} was proved in~\cite{kaneko}.
The case $n=2$ was established in~\cite{uludag}. 
In this case, $D_2^{(b)}$ coincides with $D_{n,1}^{(b)}$, which 
is a curve of genus  $\frac{1}{2}(b^2-3b+2)$ with 
$3b$ cusps of type $x^2=y^b$ and no other singularities.  
Irreducibility of $D_n^{(b)}$ is proved in Proposition~\ref{irreducible}. 
The remaining assertions of Theorem~\ref{main} are proved in Theorem~\ref{main2}.
Our method  leads naturally to the definition of braid groups of the orbifolds over $\pone$,
which we dicuss in Section~4 below. These groups were already introduced in~\cite{allcock}
for some basic cases.

\medskip\noindent
\textbf{2. Orbifolds.}
Let $M$ be a connected complex manifold, $G\subset \mbox{Aut}(M)$ a properly
discontinuous subgroup and put $N:=M/G$. 
Then the projection $\phi:M\ok N$ is a branched Galois 
covering  endowing $N$ with a map $\beta_\phi:N\ok\N$ defined by 
$\b_\phi(p):=|G_q|$ where 
$q$ is a point in $\phi^{-1}(p)$ and $G_q$ is the isotropy subgroup of $G$ at $q$. 
In this setting, the pair $(N,\b_\phi)$ is said to be uniformized by $\phi:M\ok (N,\b_\phi)$.
An \textit{orbifold} is a pair $(N,\b)$ of an irreducible normal analytic space $N$ with a 
function $\b:N\ok \N$ such that  the pair $(N,\b)$ is locally finitely uniformizable. 
A covering $\phi:(N\pr,\b\pr)\ok (N,\b)$ of orbifolds is a branched 
Galois covering $N\pr\ok N$ with 
$\b\pr=(\b\comp \phi)/\b_\phi\comp \phi$.
Note that the restriction $(N\pr,1)\ok (N,\b_\phi)$ is a uniformization of $(N,\b_\phi)$.
Conversely, let $(N,\b)$ and $(N,\c)$ be two orbifolds with $\c|\b$, and let 
$\phi:(N\pr,1)\ok (N,\c)$ be a uniformization of $(N,\c)$, e.g. $\b_\phi=\c$. 
Then  $\phi:(N\pr,\b\pr)\ok (N,\b)$ is a covering, where 
$\b\pr:=\b\comp\phi/\c\comp\phi$. 
The orbifold $(N\pr,\b\pr)$ is called the 
\textit{lifting of $(N,\b)$ to the uniformization $N\pr$ of $(N,\c)$}.
\par
Let $(N,b)$ be an orbifold,  
$B_\b:={\rm supp}(\b-1)$ and let $B_1,\dots,B_n$ be the irreducible components of $B_\b$.
Then $\b$ is constant on $B_i\eksi {\rm sing}(B_\b)$; so let $b_i$ be this number.
The \textit{orbifold fundamental group} $\orbfg(N,\b)$ of $(N,\b)$ is the group defined by
$\orbfg(N,\b):=\pi_1(N\eksi B_\b)/\ssol\mu_1^{b_1},\dots, \mu_n^{b_n}\ssag$
where $\mu_i^{b_i}$ is a meridian of $B_i$ and $\ssol\ssag$ denotes the normal closure.
An orbifold $(N,\b)$ is said to be \textit{smooth} if $N$ is smooth. 
In case $(N,\b)$ is a smooth orbifold the map $\beta$ is determined by the numbers 
$b_i$; in fact  $\b(p)$ is the order of the local orbifold fundamental group at $p$. 
Since the orbifolds to be considered in this article are exclusively smooth, 
we shall adopt the convention that such orbifolds are defined to be the pairs
$(N,B)$ where $B:=b_1B_1+\dots +b_nB_n$ is a divisor with $b_i\geq 1$.
We shall also allow $b_i$ to take infinite values, meaning that the corresponding hypersurface $B_i$
is removed from the base space $N$. If $\co:=(N,B)$ is an orbifold and $C$ a hypersurface in $N$, 
then we shall use the notation $(\co,bC)$ to denote the orbifold $(N,B+bC)$.

\medskip\noindent
\textbf{3. Discriminants.}
Let $n\geq 1$ be an integer and consider the action of the symmetric group $\Sigma_n$ on 
$(\pone)^n$.  
Let  $p_i=[u_i,v_i]\in \pone$ and let $\sigma_j$ ($j\in [0,n]$) be the 
homogeneous elementary symmetric polynomial
$$
\sigma_j(p_1,\dots,p_n):=\sum_{A\subset [1,n],\, |A|=j}
\left(
\prod_{\alpha\in A}x_\alpha
\prod_{\beta\in [1,n]\eksi A}y_\beta
\right)
$$ 
It is well known that the map  $\phi:(\pone)^n\ok \pn$  given by  
$$
\phi: (p_1,\dots, p_n):=
[\sigma_0(p_1,\dots,p_n):\dots:\sigma_n(p_1,\dots, p_n)]
$$ 
is $\Sigma_n$- invariant and gives an isomorphism 
$(\pone)^n/\Sigma_n\simeq\pn$.
\par
Let $\pi_i:(\pone)^n\ok \pone$ be the $ith$ projection map, 
$q$ a point in $\pone$, and put $F^i_q:=\pi_i\1(q)$. 
Let $\tau_{ij}\in\Sigma_n$ be the transposition exchanging 
the $i$th and $j$th coordinates of $(p_1,\dots,p_n)\in (\pone)^n$.
Since $\tau_{1i}F^1_q= F^i_q$, the hypersurface $H_q:=\phi_n(F^i_q)$ does not depend on $i$.
\begin{lemma}\label{hypps}
For any $q\in \pone$,  the hypersurface $H_q$ is a hyperplane in $\pn$.
For any set $\{q_0,\dots,q_m\}\subset \pone$ of distinct points, 
the hyperplanes $H_{q_0},\dots,H_{q_m}$ 
are in general position.
\end{lemma} 
\noindent \textit{Proof.}
Suppose without loss of generality that $i=1$. Then $H_q$ is parametrized as 
$H_q=[X_0:X_1:\dots:X_n]\in\pn$,  where $X_j=\sigma_j(q,p_2,\dots,p_n)$ and 
$ p_i\in\pone$ ($i\in [2,n]$). 
If $q=[u_1:v_1]=[x:y]$ and $p_i=[u_i:v_i]$  ($i\in [2,n]$) then one has the identity
\begin{equation}\label{discrim}
P(A,B):=\sum_{j\in[0,n]} (-1)^{n-j}\sigma_j(q,p_2,\dots,p_n)A^jB^{n-j}=
\prod_{i\in[1,n]}(u_iA-v_iB)
\end{equation}
Substitute $[A:B]=[y:x]$ in~(\ref{discrim}). 
Since the right-hand side of (\ref{discrim}) vanish at the point $(q,p_2,\dots,p_n)$, 
so does the middle term, and thus $H_q$ satisfies the linear equation
\begin{equation}\label{linear}
\sum_{j\in [0,n]} (-1)^{n-j}y^jx^{n-j}X_j =0
\end{equation}
Let $\{q_i=[x_i:y_i]\,:\,i\in[0,n]\}$ be a set of $n+1$ points.
Since the determinant of the projective Vandermonde matrix $\cv(q_0,\dots,q_n)$ given by 
$$
\cv_{i,j}(q_0,\dots,q_n):=(-1)^{n-j}y_i^jx_i^{n-j} \quad i,j\in[0,n] 
$$
vanish if and only if $q_i=q_j$ for some $i,j\in[0,n]$, the hyperplanes 
$H_{q_0},\dots,H_{q_n}$ are always in general position.
\hfill$\Box$

\medskip\noindent
The hypersurface
$$
\Delta_n:=
\{(p_1,\dots,p_n)\in (\pone)^n\,:\,p_i= p_j \,\, \text{for some}\,\, 1\leq i\neq j\leq n\}
$$ 
of $(\pone)^n$ consists of points fixed by an element of $\Sigma_n$, so that the
covering $\phi$ is branched along the hypersurface 
$D_n:=\phi(\Delta_n)$, which  is called the \textit{discriminant hypersurface}
since it is defined by the discriminant of the homogeneous polynomial $P(A,B)$.

\par
In the orbifold terminology, one has an orbifold covering
\begin{equation}\label{covering1}
\left ((\pone)^n,a\Delta_n\right)\ok (\pn, 2aD_n)
\end{equation}
\par
Let $\{q_0,\dots, q_m\}\subset\pone$ be $m+1$ distinct points, $b_0,\dots,b_m$ numbers in
$ \N\cup\{\infty\}$ and consider the orbifold 
$$
\cf(b_0,\dots,b_m):=(\pone,b_0q_0+\dots+b_mq_m)
$$
Let $ n\geq 1$ be  an integer and consider
the orbifold $\cf(b_0,\dots,b_m)^n$. Let $\cg_n$ be the orbifold 
$$
\cg_n(a,b_0,\dots,b_m):=\bigl(\cf(b_0,\dots,b_m)^n,a\Delta_n\bigr)
$$
and define the orbifold $\ch_n(a,b_0,\dots,b_m)$ as
$$
\ch_n(a,b_0,\dots,b_m):=(\pn,aD_n+b_0H_{q_0}+\dots+b_mH_{q_m})
$$
By the covering in (\ref{covering1}) and Lemma~\ref{hypps} one has the following fact 
\begin{lemma} \label{ocovering1}
There is an orbifold covering of degree $n!$
$$
\phi: \cg_n(a,b_0,\dots,b_m)\to \ch_n(2a,b_0,\dots,b_m)
$$
\end{lemma}
\noindent
In particular, for $a=1$ one has the orbifold covering
$$
\phi:\cf(\bdots)^n\simeq \cg_n(1,\bdots)\ok\ch_n(2,\bdots)
$$
The following facts are well known (see~\cite{namba}):
\begin{theorem}\label{buni}[Bundgaard-Nielsen,Fox] 
The orbifold $\cf(\bdots)$ admits a finite uniformization if $n>1$, 
$b_i<\infty$ ($1\leq i\leq m$) and if $n=2$, then $b:=b_0=b_1$. 
Let $R\ok \cf(\bdots)$ be a finite uniformization.\\
(i) $R\simeq \pone$ if $n=2, b_0=b_1<\infty$ or $n=3$, $b_0\1+b_1\1+b_2\1>1$. 
In this case , $\pone$ is also the universal uniformization. 
The groups $\orbfg\bigl(\cf(b,b)\bigr)$ and 
$\orbfg\bigl(\cf(b_0,b_1,b_2)\bigr)$ are finite
of orders $b$ and $2\bigl[b_0\1+b_1\1+b_2\1-1\bigr]\1$ respectively.\\
(ii) $R$ is of genus 1 if $n=3$, $b_0\1+b_1\1+b_2\1=1$ or $n=4$, $b_0=b_1=b_2=b_3=2$. 
Hence, $\C$ is the universal uniformization of these orbifolds.
Moreover, $\cf(\infty,\infty)$ and $\cf(2,2,\infty)$ are uniformized by $\C$. 
The corresponding orbifold fundamental groups are infinite solvable.\\
(iii) $R$ is of genus $>1$ otherwise, and the universal uniformization is 
$(\B_1)^n$, where $\B_1$ 
is the unit disc in $\C$. The corresponding orbifold fundamental groups are big 
(i.e. they contain non-abelian free subgroups).
\end{theorem}  

\noindent 
In virtue of the covering $\phi:\cf(\bdots)^n\ok\ch_n(2,\bdots)$ one has the 
following corollary.
\begin{corollary}\label{cor1}
Let $n>1$, $b_i<\infty$ ($1\leq i\leq m$) and if $n=2$, then $b_0=b_1$. 
Then the orbifold $\ch_n(2,\bdots)$ admits a finite uniformization by $R^n$, where 
$R$ is the uniformization of $\cf(\bdots)$ given in Theorem~\ref{buni}.
The orbifolds $\ch(2,\infty,\infty)$ and $\ch(2,2,2,\infty)$ are uniformized by $\C^n$.
Moreover, $\orbfg\bigl(\ch(b,b)\bigr)$ is a finite group of order $n!b^n$ and 
$\orbfg\bigl(\ch(b_0,b_1,b_2)\bigr)$ is a finite group 
of order $n!2^n\bigl[b_0+b_1+b_2-1\bigr]^{-n}$ if $b_0\1+b_1\1+b_2\1>1$.
\end{corollary}

\medskip\noindent
\textbf{4. Braid groups.}
Following and generalizing~\cite{allcock}, let us call the groups 
$$
\pure_n(a,\bdots):=\orbfg\bigl(\cg_n(a,b_0,\dots,b_m)\bigr)
$$ 
the  \textit{pure braid groups of $\cf(\bdots)$ on $n$ strands},
and the groups 
$$
\braid_n(a,\bdots):=\orbfg\bigl(\ch_n(a,\bdots)\bigr)
$$ 
the \textit{braid groups of $\cf(\bdots)$ on $n$ strands}. 
Obviously, the group $\braid_n(\abdots)$ is a quotient of $\braid_n(a\pr,b_0\pr,\dots,b_m\pr)$
provided $a|a\pr$ and $b_i|b_i\pr$ for $0\leq i\leq n$. 
The group $\braid_n(\abdots)$ is a subgroup of 
$\braid_{n+k}(a,b_{k},\dots,b_m)$ in case the equality  $a=b_0=\dots=b_{k-1}$ holds.
The group $\braid_n(2a, \bdots)$ is a normal subgroup of index $n!$ in the group
$\pure_n(a,\bdots)$. The group $\braid_n(\abdots)$ admits the presentation 
(see \cite{amram} for the case $n=2$ and \cite{paolo}\footnote{The projective relation  was kindly communicated by Paolo Bellingeri.}, \cite{paolo2}, \cite{lambro} for the general case )

\bigskip\noindent
\textit{generators} \\ $\s_1,\dots,\s_{n-1}$, $\t_0,\dots,\t_m$\\
\textit{braid relations} \\ 
$[\s_i,\s_j]= 1$, $|i-j|>1$\\ 
$\s_i\s_{i+1}\s_i=\s_{i+1}\s_i\s_{i+1}$, $1\leq i\leq n-1$\\
\textit{mixed relations}\\
$(\s_1\t_i)^2 = (\t_i \s_1)^2$,  $1 \leq i \leq m$, \\
$[\t_i,\s_j]=1$, $j\neq 1$,  $1\leq i\leq m$\\
$[\s_1 \t_i \s_1\1, \t_j] = 1$, $1 \leq i<j\leq m$, \\
\textit{projective relation }$\!\! ^1$  \\
$\s_1\s_2\dots \s_{n-1}\t_0 \cdots \t_m \s_{n-1} \dots \s_2\s_1= 1$\\
\textit{orbifold relations}\\
$\t_0^{b_0}=\dots=\t_m^{b_m}=\s_1^a=1$

\bigskip
In particular, the group $\braid_{n}(\infty,\infty)$ is the usual braid group of $\C$ 
introduced by Artin~\cite{artin}. 
The group $\braid_{n}(\infty)\simeq \braid_n(\infty,1)$ 
is the \textit{braid group of the sphere}, see~\cite{zariski}.
On the other hand, one has 
$$
\braid_1(\bdots)\simeq \langle \t_0, {\dots}, \t_m\,|\, \t_0^{b_0}= {\dots} =\t_m^{b_m}
=\t_0\cdots \t_m=1\rangle
$$
\par
In case $n=2$, the discriminant hypersurface $D_2^{(1)}$ is a smooth quadric,  
and the lines $H_{q_i}$ are tangent  to $D_2^{(1)}$ 
(see~\cite{uludag}). 
In particular, the groups $\braid_2(a,b)$ are  abelian. 
The group $\braid_2(a,b,c)$ admits the presentation
\begin{equation}\label{solvableq}
\braid_2(a,b,c)\simeq\langle \t,\s\,|\, (\t\s)^2=(\s\t)^2,\quad
 \t^{b}=(\t\s^2)^{c}=\s^a=1\rangle
\end{equation}
\begin{proposition}\label{solvable}
For $b,c<\infty$, the group $\braid_2(a,b,c)$ is a finite central extension of the triangle group 
$T_{2,a,d}:\langle \t,\s\,|\, (\t\s)^2= \t^{d}=\s^a=1\rangle$, where $d:={\rm gcd}(b,c)$. Hence,
 $\braid_2(a,b,c)$ is finite $1/a+1/b>1/2$, infinite almost solvable if $1/d+1/a=1/2$,
and big otherwise (i.e. it contains non-abelian free subgroups). The group $\braid_2(a,b,b)$ 
is of order $2b\bigl[{a}\1+{b}\1-{2}\1\bigr]\1$ if $1/a+1/b>1/2$.
\end{proposition}
\noindent\textit{Proof.}
Note that $\delta:=(\t\s)^2$ is central in $\braid_2(a,b,c)$, so that 
$(\t\s^2)^c=1$ $\Leftrightarrow $ $(\s\t\s)^c=1$ 
$\Leftrightarrow$ $ (\t\1\delta)^c=\t^{-c}\delta^c=1$. The element $\delta$ is of finite order.
Adding the relation $\delta=1$ to the presentation (\ref{solvableq}) yields the 
triangle group $T_{2,a,d}$, which is finite if $1/a+1/d>1/2$, infinite solvable if $1/a+1/d=1/2$, 
and big otherwise. In case $c=b$, one has $d=b$ and the triangle group is of order 
$2\bigl[{a}\1+{b}\1-{2}\1\bigr]\1$ if $1/a+1/b>1/2$, which shows that $\braid(a,b,b)$ is of order
$2b\bigl[{a}\1+{b}\1-{2}\1\bigr]\1$.
\hfill \qed

\medskip
Let   $R^n$ be a uniformization of the orbifold $\bigl(\cf(\bdots)\bigr)^n$. 
If $k\geq m$ then any orbifold $\cg_n(2a, c_0b_0,\dots c_mb_m,c_{m+1},\dots,c_k)$ can be lifted to $R^n$.
In case $R\simeq \pone$ or $R\simeq \C$ one obtains some arrangements associated to 
reflection groups as follows. Suppose that $q_0=[0:1]$ and $q_1=[1:0]$. 
Lifting $\cg_n(2a,cb,\infty)$ to the uniformization of $\cg_n(2,b,b)$ yields 
the orbifold $(\C^n,a\Delta^{(b)}_n+cF)$ where 
$F:=\{(X_1,\dots,X_n)\in \C^n\,:\,X_1\cdots X_n=0\}$ and $\Delta^{(b)}_n$ 
is the lifting of the superdiagonal
$$
\Delta^{(b)}_n
:=\{(X_1,\dots,X_n)\in \C^n\,:\,\psi_b(p_i)= \psi_b(p_j) 
\,\, \text{for some}\,\, 1\leq i\neq j\leq n\}
$$
with $\psi_b(X)=X^b$ if $b<\infty$ and $\psi_\infty(X)=\exp(2\pi {\rm i }X)$.
Setting $b=2$ in this construction identifies the group 
$\braid_n(\infty,\infty,\infty)$ with the Artin group corresponding 
to the diagrams $B_n$ (see~\cite{allcock}).

The groups $\braid_{2}(a,b,c,d)$ admits the simplified presentation (see~\cite{uludag})
$$
\braid_2(a,b,c,d) \simeq 
\left\langle \t,\r,\s\,\left|\,
                        \begin{array}{l}
                         (\t\s)^2=(\s\t)^2,\; (\r\s)^2=(\s\r)^2,\;[\r,\t]=1,\\
                        \t^{b}=(\s\t\s\r)^{d}=\r^{c}=\s^a=1
                        \end{array}
                        \right. \right\rangle
$$
We summarized the known information about the orbifolds $\ch$ and the 
corresponding braid groups in Table 1 below. Suppose that if $(n,m)=(2,1)$ then 
$b_0=b_1$. We believe that the group
$\braid_n(\abdots)$ is finite if 
$$
\frac{2(n-1)}{a}+\sum_{i\in[0,m]}\frac{1}{b_i}>n+m-2,
$$
infinite solvable if the equality holds, and big otherwise.

\bigskip
{\tiny \begin{tabular}{|l|c|l|l|}
\hline
\textbf{Orbifold} &\textbf{Uniformization} & \textbf{Braid group} & \textbf{Reference}\\
\hline
${\ch_n(2)}$& $(\pone)^n$        & $n!$ & Cor.~\ref{cor1}\\
\hline
${\ch_n(2,b,b)}$& $(\pone)^n$& $n!b^n$&Cor.~\ref{cor1}\\
\hline
\!\!\!\!$\col{\ch_n(2,b,c,d)}{\bigl(1/b+1/c+1/d>1\bigr)}$ 
& $(\pone)^n$&${n!2^n}{\bigl[\frac{1}{b}+\frac{1}{c}+\frac{1}{d}-1\bigr]^{-n}}$ & Cor.~\ref{cor1}\\
\hline
\!\!\!\!$\col{\ch_n(2,b,c,d)}{\bigl(1/b+1/c+1/d=1\bigr)}$
&$\C^n$&{Crystallographic}& Cor.~\ref{cor1}\\
\hline
$\ch_n(2,2,2,2,2)$&$\C^n$&Crystallographic&Cor.~\ref{cor1}\\
\hline
\!\!\!\!$\col{\ch_n(2,\bdots)}{{\rm (otherwise)}}$& $(\B_1)^n$ & Linear & Cor.~\ref{cor1}\\
\hline
$\ch_n(\infty,\infty,\infty)$& -& $B_n $-Artinian &~\cite{brieskorn}\\
\hline
$\ch_2(\infty,\infty,\infty,\infty)$& -& $\widetilde{C}_2$-Artinian&~\cite{brieskorn}\\
\hline
\!\!\!\!$\col{\ch_2(a,b,b)}{(1/a+1/b>1/2)}$& - &$2b\bigl[\frac{1}{a}+\frac{1}{b}-\frac{1}{2}\bigr]\1$ &Prop.~\ref{solvable}\\
\hline
\!\!\!\!$\col{\ch_2(a,b,b)}{(1/a+1/b=1/2)}$& - &$\infty$ almost solvable &Prop.~\ref{solvable},
\cite{uludag}\\
\hline 
${\ch_3(a,\infty)}\, {(a=3,4,5)}$& - & 24, 96, 600 & \cite{coxeter}\\
\hline
${\ch_n(3,\infty)}\, {(n=4,5)}$ &-& 648, 155520 & \cite{coxeter}\\
\hline
$\ch_3(\infty,2)$&-&192&Maple\\
\hline
${\ch_4(a)}\,{(a=4,5)}$&-&192, 60&Maple\\
\hline
$\ch_5(4)$&-&120&Maple\\
\hline
$\ch_2(a,2,2,2)$ &K3 $(a=4)$ &$ 4a^3$ & \cite{uludag}\\
\hline
$\ch_2(3,3,2,2)$ &K3 & 576 & \cite{uludag}\\
\hline
\!\!\!\!$\begin{array}{ll}
        \ch_2(3,3,4,4)&\ch_2(4,4,4,4)\\
        \ch_2(3,6,6,2)&\ch_2(3,3,3,6)\\
\end{array} $
&$\B_2$& Picard Modular & \cite{holzapfel},\cite{uludag}\\
\hline
$\ch_2(3,3,4,2)$&$\B_1\times\B_1 ?$-& Unknown & \cite{uludag}\\
\hline
\end{tabular}}

\medskip
\begin{center}
{\bf Table 1.}
\end{center}

\bigskip\noindent
\textbf{5. Another covering of $\ch(\abdots)$.}
Let $b\in \N$ be an integer and 
consider the orbifold $\ck_n(b):=(\pn,bH_{q_0}+\dots+bH_{q_n})$. 
By Lemma~\ref{hypps}, the hyperplanes $H_{q_0},\dots, H_{q_n}$ are in general position. 
It is well known that the universal uniformization of this orbifold is $\pn$. Applying a projective 
transformation one may assume that the hyperplanes
$H_{q_i}$ are given by the equations $Y_i=0$ where $[Y_0:\dots:Y_n]\in \pn$. 
In this case the uniformization
$\psi_b: \pn\ok\ck_n(b)$ is nothing but the map 
$$
[Y_0:\dots:Y_n]=\psi_b([Z_1:\dots:Z_n])=[Z_1^b:\dots:Z_n^b]
$$
It is clear that the orbifold $\ch_n(a,bb_0,\dots,bb_n)$ lifts to the uniformization of 
$\ck_n(b)$. Put $D_n^{(b)}:=\psi_b\1(D_n)$,  denote $M_{q_i}:=\psi\1(H_{q_i})$ 
and define the orbifold
$$
\cl_n^{(b)}(\abdots):=(\pn, aD_n^{(b)}+b_0M_{q_0}+\dots b_nM_{q_n})
$$
to be this lifting. In case $n=2$ these liftings were studied in~\cite{uludag}. 
For $n>2$ the following proposition is valid:
\begin{proposition}
For $n>2$ and $b\geq 2$  
the orbifolds $\cl_n^{(b)}(2,\bdots)$ are uniformized by $(\B_1)^n$ except the 
orbifold $\cl_3^{(2)}(2)$, which is uniformized by $\C^3$. 
\end{proposition}
\noindent\textit{Proof.} There is an orbifold covering
$$
\cl_n^{(b)}(2,\bdots)\longrightarrow \ch_n(a,bb_0,\dots,bb_n)
$$
The claim follows, since by Corollary~\ref{cor1} 
the latter orbifold is uniformized by $\C^3$ if $b=2$, $n=3$, $b_0=\dots=b_n=1$
and by $(\B_1)^n$ otherwise. \hfill $\Box$

\medskip\noindent
For $k\in[1,n]$, define the $k$-dimensional subvarietiy $\Delta_{n,k}$ of $\Delta_{n}$ by 
$$
\Delta_{n,k}:=\{(p_1,p_2,\dots,p_n)\in (\pone)^n\,:
\,p_{k}=p_{k+1}=\dots=p_n\}\simeq (\pone)^k
$$
Thus, $\Delta_{n,n-1}$ is an irreducible component of $\Delta_{n}$ and  
$\Delta_{n,1}$ is the diagonal in $(\pone)^n$.  
The subgroup of $\Sigma_n$ acting on $\Delta_{n,k}$ is a symmetric group $\Sigma_{k-1}$,
so that $D_{n,k}:=\pone\times{\mathbb P}^{k-1}$. These varieties admits the 
parametrizations
\begin{equation}\label{param}
D_{n,k}:[X_0:\dots:X_n]\in\pn\quad X_j=\s_j(p_1,\dots,p_n),\quad p_k=\dots=p_n
\end{equation}
In particular, the curve  $D_{n,1}$ is a rational normal curve parametrized as
$$
\left[\comb{n}{0}v^n,\comb{n}{1}uv^{n-1},\dots,\comb{n}{n}u^n\right]\quad ([u:v]\in \pone)
$$ 
Applying the projective transformation $\cv(q_0,\dots,q_n)$ 
to the parametrizations (\ref{param}) gives the parametrization
$D_{n,k}:[Y_0:\dots:Y_n]\in\pn$, where 
\begin{equation}\label{param1}
\sum_{j\in [0,n]} (-1)^{n-j}y_i^jx_i^{n-j}\sigma_j(p_1,\dots,p_n), \quad p_k=\dots=p_n 
\end{equation}
Let $p_i=[u_i:v_i]$ and let $[u:v]=[u_k:v_k]=\dots=[u_n:v_n]$. 
In virtue of the identity (\ref{discrim})  one has the parametrizations
$D_{n,k}:[Y_0:\dots:Y_n]\in\pn$ where
\begin{equation}\label{param2}
 Y_j=(uy_j-vx_j)^{n-k+1}\prod_{i\in[1,k-1]}(u_iy_j-v_ix_j)
\end{equation}
In particular, the curve $D_{n,1}$ is parametrized as
\begin{equation}\label{normal}
D_{n,1}:\bigl[(uy_0-vx_0)^n:\dots:(uy_n-vx_n)^n\bigr]
\end{equation}
The varieties $D_{n,k}^{(b)}$ are parametrized as
\begin{equation}\label{param3}
D_{n,k}^{(b)}:[Z_0:\dots:Z_n]\quad Z_j^b=
(uy_j-vx_j)^{n-k+1}\prod_{i\in [1,k-1]}(u_iy_j-v_ix_j)
\end{equation}
Note that the parametrizations (\ref{param2}) and (\ref{param3}) are not generically 
one-to-one unless $k\leq 2$, since  (\ref{param2}) is a map $(\pone)^k\ok D_{n,k}$. 
\begin{proposition}\label{irreducible}
(i) The curve  $D_{n,1}^{(b)}$ is irreducible if and only if ${\rm gcd}(n,b)=1$. 
Hence, the subvarieties $D_{n,k}^{(b)}$ are irreducible if ${\rm gcd}(n,b)=1$.  
\end{proposition}

\noindent\textbf{Definition.}
Let $t\in \Z$ and $\psi_t$ be the map
$$
\psi_t: [Z_0:\dots:Z_n]\in\pn\ok [Z_0^t:\dots Z_n^t]\in \pn 
$$
Let $V\subset \pn$ be a subvariety and $r,s\in\Z$ such that $s>1$. 
Then $V^{(r/s)}$ is the subvariety of $\pn$ defined as 
$$
V^{(r/s)}:=(\psi_r\1\comp\psi_s)(V)
$$
In particular, $V^{(r/r)}$ is the orbit of $V$ under the $(\Z/(r))^n$-action on $\pn$. 

\medskip\noindent\textit{Proof of the Proposition.} 
The parametrization~(\ref{normal}) shows that $D_{n,1}\simeq L^{1/n}$, 
where $L$ is a line $\subset\pn$ in general position 
with respect to $\psi_n$, in other words $L$ intersects the hyperplane arrangement 
$Z_0\dots Z_n=0$ transversally at smooth points. 
Hence there is a surjection of fundamental groups
\begin{equation}\label{surjection}
\pi_1(L\eksi \{\tilde{q}_0,\dots,\tilde{q}_n\})\onto \pi_1(\pn\eksi \{Z_0,\dots,Z_n\})
\end{equation}
where $\tilde{q}_i:=Z_i\cap L$. Let $\cm(b)$, $\ck(b)$ be the orbifolds
$$
\cm(b):=(L,b\tilde{q}_0+\dots+b\tilde{q}_n),\quad
\ck(b):=(\pn, bZ_0+\dots+bZ_n)
$$
Then (\ref{surjection}) induce a surjection of orbifold fundamental groups
$$
\orbfg\bigl(\cm(b)\bigr)\onto \orbfg\bigl(\ck(b)\bigr)
$$
(one may say: $\cm(b)$ is a sub-orbifold of $\ck(b)$). 
This shows that the curve $L^{(b)}$ is irreducible and is a uniformization of $\cm(b)$.
Since ${\rm gcd}(n,b)=1$, one has $D_{n,1}^{(b)}=L^{(b/n)}$, 
showing that $D_{n,1}^{(b)}$ is irreducible. Note that  $D_{n,1}^{(b)}$ is the maximal abelian 
orbifold covering of $\cm(b)$. Irreducibility of $D_{n,k}^{(b)}$
follows since $D_{n,1}^{(b)}$ is a subvariety of $D_{n,k}^{(b)}$. \hfill\qed

\medskip\noindent
Let $\co(b)$ be the orbifold 
$$
\co(b):=(D_{n,1},b\bar{q}_0+\dots+b\bar{q}_n),
$$
where $\bar{q}_i:=Y_i\cap D_{n,1}$. 
The orbifold $\co(b)$ is identified via the covering $\phi $ with the orbifold
$$
\cp(b):=(\Delta_{n,1},bq_0\pr+\dots+bq_n\pr),
$$
where this time $q_i\pr:=\phi \1(\bar{q_i})$. 
In turn, $\co(b)$ is identified with the orbifold $\cf(b,\dots,b)$ via the coordinate projection.  
By the proof of  Proposition~\ref{irreducible}, these orbifolds are identified with the 
orbifold $\cm(b)$ in case $(n,b)=1$. 

\begin{theorem}\label{main2}
Let ${\rm gcd}(n,b)=1$. Then there is a finite uniformization
$\xi_n:(D_{n,1}^{(b)})^n\longrightarrow \cl_n^{(b)}(2)$ which is of degree $n!b^{n^2-n}$.
\end{theorem}

\noindent\textit{Proof.}
One has the diagram
$$
\begin{array}{ccc}
&\xi_n&\\
\cl_n^{(b)}(2)&{\longleftarrow\!\!\!-\!\!\!-\!\!\!-}&(D_{n,1}^{(b)})^n\\
&&\\
\psi_b\,\,\,\Bigg\downarrow&&\Bigg\downarrow\,\,\,\zeta_b\\
&\phi_n&\\
\ch_n(2,b,\dots,b)&{\longleftarrow\!\!\!-\!\!\!- \!\!\!-}&\co(b)^n
\end{array}
$$

\medskip\noindent
where $\zeta_b:(D_{n,1}^{(b)})^n\ok \co(b)^n$ is the maximal abelian orbifold covering
and $\xi_n$ is to be shown to be a branched Galois covering of degree $n!b^{n^2-n}$.
It suffices to show that the group $H:=(\phi_n\comp \zeta_b)_* \pi_1\bigl((D_{n,1}^{(b)})^n\bigr)$ 
is a normal subgroup of  $K:=(\psi_b)_*\orbfg\bigl(\cl_n^{(b)}(2)\bigr)$. 
Let $\sigma$ be a meridian of $D_n$. Then since
$\orbfg\bigl(\ch_n(2,b,\dots,b)\bigr)/\ssol \sigma\ssag\simeq
\orbfg\bigl(\ck_n(b)\bigr)\simeq \bigl(\Z/(b)\bigr)^n$
is the Galois group of $\psi_b$, the group $K$ is the normal subgroup of 
$\orbfg\bigl(\ch_n(2,b,\dots,b)\bigr)$ generated by $\sigma$, i.e. $K\simeq \ssol \sigma\ssag$.
The group $\orbfg\bigl(\ch_n(2,b,\dots,b)\bigr)\bigr)/K$ being abelian, one has $[\t_i,\t_j] \in K$ 
for $i,j\in[0,n]$. On the other hand one has
$$
\orbfg\bigl(\ch_n(2,b,\dots,b)\bigr)/\ssol\t_0,\dots,\t_n\ssag\simeq 
\orbfg\bigl(\ch_n(2)\bigr)\simeq \Sigma_n
$$ 
Since $\Sigma_n$ is the Galois group of $\phi_n$, one has
$\phi_*\co(b)^n\simeq\ssol\t_0,\dots,\t_n\ssag.$
Since $\zeta_n$ is the maximal abelian orbifold covering,
one has $H\simeq \ssol [\t_i,\t_j]\ssag$. 
This shows that $H$ is a normal subgroup of $K$. Since $\deg(\zeta_b)=b^{n^2} $, 
$\deg(\Phi_n)=n!$ and $\deg(\psi_b)=b^n$, one has 
$$
\deg(\xi_n)=\frac{\deg(\zeta_b) \deg(\phi_n)}{\deg(\psi_b)}=n!b^{n^2-n} 
$$
The euler number of $D_{n,1}^{(b)}$ is easily computed by Riemann-Hurwitz formula.
\hfill$\Box$

\bigskip\noindent\textbf{6. Remarks.}
Consider the restriction of $D_{n,k}$ to the $n-k+1$ dimensional linear subspace 
$M_{n-k+1}:=\{ [Y_0:\dots : Y_n]\in \pn \, |  \, Y_{n-k+2}=\dots=Y_n=0\}$ of $\pn$. 
Setting  $[u:v]=[x_{n}:y_{n}]$ and $[u_i:v_i]=[x_{n-i}:y_{n-i}]$ for $i\in[1,k-2]$ in~(\ref{param2})
we see that $D_{n,k}$ has a 1-dimensional linear component $L$ 
in $M_{n-k+1}\simeq {\mathbb P}^{n-k+1}$, parametrized as
$[Y_0:\dots :Y_{n-k+1}] \in M_{n-k+1}$ where
$$
 Y_l=(u_{n}y_l-v_{n}x_l) (x_ny_l-y_nx_l)^{n-k+1}\prod_{i\in[2,k-1]}(x_{n-i}y_l-y_{n-i}x_l)
$$
for $l\in[0,n-k+1]$ and $[u_n:v_n]\in\pone$. 
It is readily seen that there are $k-1$ such lines. In case $[u_i:v_i]=0$ for 
$i\in [1,k-1]$, one has the curve $C$ in $D_{n,k}\cap M_{n-k+1}$ parametrized as
$[Y_0:\dots :Y_{n-k+1}] \in M_{n-k+1}$ where 
$$
 Y_l=(uy_l-vx_l)^{n-k+1}\prod_{i\in[1,k-1]}(x_{n-i}y_l-y_{n-i}x_l)
$$
for $l\in[0,n-k+1]$ and $[u:v]\in\pone$, which shows that $C$ is the curve 
$E^{(1/n-k+1)}$ for some line $E$ in ${\mathbb P}^{n-k+1}$. 
The lines $L$ are tangent to $C$ with multiplicity $n-k+1$. 
In case $k=n-1$, one has $M_{n-k+1}\simeq\ptwo$, and one obtains an arrangement 
of a quadric $C$ with $n-2$ tangent lines. The lines $Y_0=0$, $Y_1=0$ and $Y_2=0$ 
are also tangent to this quadric. 
\par
From these considerations it is easy  to obtain a description of the intersection 
of $D_{n,k}^{(b)}$ with ${\mathbb P}^{n-k+1}\simeq Z_{n-k+2}=\dots=Z_n=0$. 
For $D_{3,2}^{(2)}$, this is the arrangement of a quadric with four tangent lines.
\par
Let $H\subset \pn$ be a hyperplane. The intersection  $H^{(1/2)}\cap M_2$ is a quadric,
tangent to the lines $Y_0=0$, $Y_1=0$ and $Y_2=0$, which is very similar 
to the intersections $D_n\cap M_2$.  In contrast with this, there is the following fact:
In a recent article~\cite{katz}, it was proved that 
the dual of $D_n$ is one dimensional (we believe that $D_{n,k}$ and $D_{n,n-k}$ are duals),
whereas it is easy  to show that $H^{(r/s)}$ and $H^{(r/r-s)}$ are duals, so that the dual 
of $H^{(1/2)}$ is the cubic hypersurface $H^{(-1)}$. Note also that $D_n$ is of degree $2(n-1)$,
whereas $H^{(1/2)}$ is of degree $2^{n-1}$.
It is of interest to know more about the varieties $D_{n,k}^{(r/s)}$ and their duals.

\bigskip\noindent\textbf{Appendix: The curves $L^{(r/s)}$.}
In $\ptwo$, many interesting curves appears as $L^{r/s}$.
For example, $L^{1/2}$ is the curve $D_{2,1}$, 
a quadric tangent to the coordinate lines, $L^{3/2}\simeq D_{2,1}^{3}$ is a nine cuspidal
sextic, $L^{2/3}$ is a Zariski sextic with 4 nodes and 6 cusps, $L^{-1/2}\simeq D_{2,1}^{-1}$ 
is a three cuspidal quartic, $L^{(-1)}$ is a quadric passing through the intersection points 
of the coordinate lines. 
\begin{proposition}
If $r,s\geq 0$ are coprime integers , 
then $L^{r/s}$ is an irreducible curve of degree $sr$ 
and genus $(r-1)(r-2)/2$, 
with $3r$ points of type $x^r=y^s$ and $r^2(s-1)(s-2)/2$ nodes. 
\end{proposition} 

\noindent\textit{Proof.} 
We begin by proving that the curves $L^{1/s}$ are nodal.
For this, it suffices
to show that the orbit of $L$ under the action of the group 
$\Z/(s)\oplus\Z/(s)$ has only double points on $\ptwo\eksi\{xyz=0\}$. 
If $\omega:=e^{2\pi i/s}$, then the orbit of $L$ consists of the lines 
$L_{ij}:=a\omega^ix+b\omega^jy+cz=0$ for $1\leq i,j\leq s$. 
Suppose that  no pairs of 
lines among the lines $L_{i,j}$, $L_{k,l}$, $L_{p,q}$ meet on $xyz=0$. 
Then they meet at a point $\notin \{xyz=0\} $ only if the determinant of the 
matrix 

\medskip
\begin{center}
\begin{tabular}{|ccc|}
$a\omega^i$&$ b\omega^j$& $c$\\
$a\omega^k$& $b\omega^l$& $c$\\
$a\omega^p$& $b\omega^q$& $c$
\end{tabular}
\end{center}

\noindent\medskip
vanish. Since $abc\neq 0$, this is equivalent to the vanishing of 

\medskip
\begin{center}
det
\begin{tabular}{|cc|}
$\omega^\alpha-1$&$\omega^\beta-1$\\
$\omega^\gamma-1$&$\omega^\theta-1$
\end{tabular}
\end{center}

\medskip\noindent
where $\alpha:=k-i$, $\beta:=l-j$, $\gamma:=p-i$ and $\theta:=q-j$.
The integers $\alpha$, $\beta$, $\gamma$, $\theta$ are 
not multiples of $s$ by hypothesis. 
Then vanishing of the determinant implies
$$
\frac{(\omega^\alpha-1)(\omega^\theta-1)}{(\omega^\beta-1)(\omega^\gamma-1)}=1
\ise
\frac{(\omega^{\alpha/2}-\omega^{-\alpha/2})(\omega^{\theta/2}-\omega^{-\theta/2})}{(\omega^{\beta/2}-\omega^{-\beta/2})(\omega^{\gamma/2}-\omega^{-\gamma/2})}=
\omega^{(\beta+\gamma-\alpha-\theta)/2}
$$
Since the left-hand side of the latter expression is real, so must be the 
right-hand side. Therefore
$$
\mbox{Im}(e^{\pi i (\beta+\gamma-\alpha-\theta)/s})=0\ise s|\beta+\gamma-\alpha-\theta.
$$
But this means that there is a pair of lines meeting at $z=0$, contradiction. 
This shows that the curves $L^{1/s}$ are nodal.
\par
Since $L^{1/s}$ is a rational curve of degree $s$, 
it must have $(s-1)(s-2)/2$ nodes.
Since $L^{r/s}=\phi_r\1(L^{1/s})$, the number of nodes of $L^{r/s}$ is 
$r^2(s-1)(s-2)/2$. Obviously, three flex points of $L^{1/s}$ are lifted as 
$3r$ cusps of type $x^r=y^s$. The genus of $L^{r/s}$ can be calculated by the
genus formula, or by noting that the curves $L^{r/s}$ are 
coverings of $L^{1/s}$ branched at these three flex points,
with the branching index $r$. $\Box$

\noindent
{\sc Galatasaray University,
Department of Mathematics,
Ortak{\"o}y/{\.I}stanbul,
Turkey}\\
e-mail adress: muludag@gsu.edu.tr
\end{document}